\documentclass[11pt,a4paper]{article}
\usepackage[T2A]{fontenc}
\usepackage[cp1251]{inputenc}
\usepackage{amssymb}
\usepackage{amsthm}
\newtheorem{Theorem}{Theorem}
\newtheorem{Proposition}{Proposition}
\newtheorem{Definition}{Definition}

\newcommand{\pr}{\mbox{pr}_1\,}
\newcommand{\spr}{\mbox{\scriptsize pr}_1\,}

\newcommand{\op}[1]{\mathop{\oplus}\limits_{\phantom{.}#1}}
\newcommand{\opp}[2]{\mathop{\oplus}\limits_{\phantom{.}#1}^{\phantom{.}#2}}

\title{Representations of $(2,n)$-semigroups by multiplace
functions}
\author{Wies{\l}aw A. Dudek and Valentin S. Trokhimenko}
\date{}

\begin{document}
\sloppy \maketitle


\begin{abstract}
We describe the representations of $(2,n)$-semigroups, i.e.
groupoids with $n$ binary associative operations, by partial
$n$-place functions and prove that any such representation is a
union of some family of representations induced by Schein's
determining pairs.
\end{abstract}

Let $A^n$ be the $n$-th Cartesian product of the set $A$. Any
partial mapping from $A^n$ into $A$ is called a {\it partial
$n$-place function}. The set of all such functions is denoted by
${\mathcal F}(A^n,A)$. The set of all {\it full $n$-place
functions} on $A$, i.e. mappings defined for every
$(a_1,\ldots,a_n)\in A^n$, is denoted by ${\mathcal T}(A^n,A)$.
Obviously ${\mathcal T}(A^n,A)\subset{\mathcal F}(A^n,A)$. Note
that in many papers full $n$-place functions are called {\it
$n$-ary operations}.

On ${\mathcal F}(A^n,A)$ we define $n$ binary compositions
$\,\op{1\;},\ldots,\op{n}\,$ putting
\begin{equation}\label{1}
(f\op{i\;}g)(a_1,\ldots,a_n)=
f(a_1,\ldots,a_{i-1},g(a_1,\ldots,a_n),a_{i+1},\ldots,a_n)
\end{equation}
for all $f,g,g_1,\ldots,g_n\in {\mathcal F}(A^n,A)$ and
$(a_1,\ldots,a_n)\in A^n$, where the left and right hand of
(\ref{1}) are defined or not defined simultaneously. Since all
compositions $\op{1\;},\ldots,\op{n}$ are associative operations,
algebras of the form $(\Phi;\op{1\;},\ldots,\op{n})$, where
$\Phi\subset\mathcal{F}(A^n,A)$, are called {\it
$(2,n)$-semigroups of $n$-place functions}. If
$\Phi\subset\mathcal{T}(A^n,A)$, then we say that
$(\Phi;\op{1\;},\ldots,\op{n})$ is a {\it $(2,n)$-semigroup of
full $n$-place functions} (or {\it $n$-ary operations}).

The study of such compositions of functions were initiated by Mann
\cite{Man} for binary operations and continuated by other authors
(cf. for example \cite{Bel}, \cite{Tro}, \cite{Zar}). Nowadays
such defined compositions are called {\it Mann's compositions} or
{\it Mann's superpositions}. Mann's compositions of $n$-ary
operations are described in \cite{Yak}. Abstract algebras
isomorphic to some sets of operations closed with respect to these
compositions are described in \cite{Sok}. The sets of partial
functions closed with respect to these compositions and some
additional ope\-rations are characterized in \cite{DT1}.

According to the general convention used in the theory of $n$-ary
systems, the sequence $\,a_i,a_{i+1},\ldots,a_j$, where
$i\leqslant j$, \ can be written in the abbreviated form as
$\,a_i^j$ \ (for \ $i>j$ \ it is the empty symbol). In this
convention (\ref{1}) can be written as
\[
(f\op{i\;}g)(a_1^n)=f(a_1^{i-1},g(a_1^n),a_{i+1}^n).
\]

The algebra $(G;\op{1\;},\ldots,\op{n})$, where
$\op{1\;},\ldots,\op{n}$ are associative binary operations on $G$,
is called a {\it $(2,n)$-semigroup}. Each its homomorphism into
some $(2,n)$-semigroup of $n$-place functions ($n$-ary operations)
is called a {\it representation by $n$-place functions}
(respectively, {\it by $n$-ary operations}). We say that a
representation is {\it faithful} if it is an isomorphism. A
$(2,n)$-semigroup for which there exists faithful representations
is called \textit{representable}.

In the sequel, all expressions of the form
$(\cdots((x\op{i_1}y_1)\op{i_2}y_2)\cdots)\op{i_k}y_k$ are denoted
by $x\op{i_1}y_1\op{i_2}\cdots\op{i_k}y_k$ or by
$x\opp{i_1}{i_k}y_1^k$. The symbol $\mu_i(\opp{i_1}{i_s}x_1^s)$
denotes an element \ $x_{i_k}\!\opp{i_{k+1}}{i_s}\!x_{k+1}^{s}$ \
if $i=i_k$\ and\ $i\neq i_p$ for all $p<k\leqslant s$. If\ $i\neq
i_p$\ for all $i_p\in\{i_1,\ldots,i_s\}$,\ this symbol is empty.
For example, $\mu_1(\op{2}x\op{1\;}y\op{3\,}z)=y\op{3\,}z$,
$\mu_2(\op{2}x\op{1\;}y\op{3\,}z)=x\op{1\;}y\op{3\,}z$,
$\mu_3(\op{2}x\op{1\;}y\op{3}z)=z$. The symbol
$\mu_4(\op{2}x\op{1\;}y\op{3\,}z)$ is empty.

For $f,g_1,\ldots,g_n\in\mathcal{F}(A^n,A)$, by $f[g_1\ldots g_n]$
(or shortly by $f[g_1^n]$) we denote the so-called {\it Menger
superposition} \cite{ST} of $n$-place functions, which is defined
by the equality
\begin{equation}\label{2}
f[g_1\ldots g_n](a_1^n)=f(g_1(a_1^n),\ldots,g_n(a_1^n)),
\end{equation}
where $a_1,\ldots,a_n\in A$. It is assumed that the left and right
hand of (\ref{2}) are defined or not defined simultaneously. By
$I^n_1,\ldots,I^n_n\,$ we denote the $n$-place {\it projectors},
i.e. $n$-place functions from $\mathcal{T}(A^n,A)$ such that
\begin{equation}\label{3}
  I^n_i(a_1,\ldots,a_n)=a_i
\end{equation}
for all $a_1,\ldots,a_n\in A$. It is not difficult to verify, that
for all $n$-place functions defined on $A$ we have
\begin{eqnarray}
&&\label{4}f[I^n_1\ldots I^n_n]=f,\\[4pt]
&&\label{5}I^n_i[g_1\ldots g_n]=g_i\circ\bigtriangleup_{\spr
g_1\cap\ldots\cap\,\spr g_n},\\[4pt]
&&\label{6}f\op{i\;}g=f[I^n_1\ldots I^n_{i-1}\,g\,I^n_{i+1}\ldots
I^n_n],\\[4pt]
&&\label{7}f[g_1^n][h_1^n]=f[g_1[h_1^n]\ldots
g_n[h_1^n]],\\[4pt]
&&\label{8}(f\op{i\;}g)[h_1^n]=f[h_1^{i-1}\,g[h_1^n]\,h_{i+1}^n],\\[4pt]
&&\label{9}f[g_1^n]\op{i\;}h=f[(g_1\op{i\;}h)\ldots
(g_n\op{i\;}h)],\\[4pt]
&&g\op{i\,}I_{i}^{n}=I_{i}^{n}\op{i\;}g=g,
\label{10} \\[4pt]
&&I_{k}^{n}\op{i\;}g=I_{k}^{n}\circ \bigtriangleup _{{\rm
{\scriptsize pr}}_{1}\,g} \label{11}
\end{eqnarray}
for all $i,k\in\{1,\ldots ,n\}$ and $k\ne i$, where
$\triangle_H=\{ (h,h)\,|\,h\in H\}$ and pr$_1\,g$ denotes the
domain of $g$.

\begin{Proposition}\label{P1}
For all $f,g_1,\ldots,g_n\in \mathcal{F}(A^n,A)$ and
$\op{i_1},\ldots,\op{i_s}$ we have
\begin{equation}\label{12}
f\opp{i_1}{i_s}g_1^s=f[(I^n_1\opp{i_1}{i_s}g_1^s)\ldots
(I^n_n\opp{i_1}{i_s}g_1^s)].
\end{equation}
\end{Proposition}
\begin{proof} We prove (\ref{12}) by induction. For $s=1$, by
(\ref{6}) and (\ref{11}), we have
\[
\begin{array}{lll}
&f\op{i_1}g_1= f[I^n_1\ldots
I^n_{i_1-1}\,g_1\,I^n_{i_1+1}\ldots I^n_n] \\[5pt]
&=f[(I^n_1\circ\bigtriangleup_{\spr g_1})\ldots
(I^n_{i_1-1}\circ\bigtriangleup_{\spr
g_1})g_1\,(I^n_{i_1+1}\circ\bigtriangleup_{\spr g_1})\ldots
(I^n_{n}\circ\bigtriangleup_{\spr g_1})] \\[5pt]
&=f[(I^n_{1}\op{i_1}g_1)\ldots
(I^n_{i_1-1}\op{i_1}g_1)(I^n_{i_1}\op{i_1}g_1)(I^n_{i_1+1}\op{i_1}g_1)\ldots
(I^n_{n}\op{i_1}g_1)] \\[5pt]
&=f[(I^n_{1}\op{i_1}g_1)\ldots (I^n_{n}\op{i_1}g_1)].
\end{array}
\]
Thus, for $s=1$, the condition (\ref{12}) is valid.

Assume that it is valid for $s=k$, i.e. that
\[
f\op{i_{1}}^{i_{k}}g_{1}^{k}=f[(I_{1}^{n}\op{i_{1}}^{i_{k}}g_{1}^{k})\ldots
(I_{n}^{n}\op{i_{1}}^{i_{k}}g_{1}^{k})].
\]
Then, according to this assumption and (\ref{9}), we obtain
\[\arraycolsep=.5mm\begin{array}{lll}
f\op{i_{1}}^{i_{k+1}}g_{1}^{k+1}&=(f\op{i_{1}}^{i_{k}}g_{1}^{k})
\!\op{i_{k+1}}\!g_{k+1}=f[(I_{1}^{n}\op{i_{1}}^{i_{k}}g_{1}^{k})\ldots
(I_{n}^{n}\op{i_{1}}^{i_{k}}g_{1}^{k})]\!\op{i_{k+1}}\!g_{k+1}
 \\[4pt]
&=f[(I_{1}^{n}\op{i_{1}}^{i_{k}}g_{1}^{k})\!\op{i_{k+1}}\!g_{k+1}\;\ldots\;
(I_{n}^{n}\op{i_{1}}^{i_{k}}g_{1}^{k})\!\op{i_{k+1}}\!g_{k+1}]
 \\[4pt]
&=f[(I_{1}^{n}\op{i_{1}}^{i_{k+1}}g_{1}^{k+1})\,\ldots\,
(I_{n}^{n}\op{i_{1}}^{i_{k+1}}g_{1}^{k+1})],
\end{array}
\]
which proves (\ref{12}) for $s=k+1$. So, on the basis of the
principle of mathematical induction, it is valid for all natural
$s$.
\end{proof}

\begin{Proposition}\label{P2}
The following implication
\begin{equation}
\bigwedge_{i=1}^{n}\left(\mu_{i}(\op{i_{1}}^{i_{s}}g_{1}^{s})=
\mu_{i}(\op{j_{1}}^{j_{k}}h_{1}^{k})\right) \longrightarrow
f\op{i_{1}}^{i_{s}}g_{1}^{s}=f\op{j_{1}}^{j_{k}}h_{1}^{k}
\label{13}
\end{equation}
is valid for all
$\,f,g_{1},\ldots,g_{s},h_{1},\ldots,h_{k}\in\mathcal{F}(A^{n},A)$
and $\,i_{1},\ldots ,i_{s},j_{1},\ldots ,j_{k}\in \{1,\ldots ,n\}$
\end{Proposition}
\begin{proof} Assume that the premise of (\ref{13}) is satisfied.
Then $\{i_{1},\ldots ,i_{s}\}=\{j_{1},\ldots ,j_{k}\}$. Indeed,
for $i\not\in\{i_{1},\ldots ,i_{s}\}$, the symbol $\mu
_{i}(\op{i_{1}}^{i_{s}}g_{1}^{s})$ is empty. So, $\mu
_{i}(\op{j_{1}}^{j_{k}}h_{1}^{k})$ also is empty. Thus $i\not\in
\{j_{1},\ldots ,j_{k}\}$.

If $\,a_{1}^{n}\in {\rm pr}_{1}\,(f\op{i_{1}}^{i_{s}}g_{1}^{s})$,
then, according to (\ref{12}), we have\vspace{-2mm}
\[\arraycolsep=.5mm\begin{array}{lll}
(f\op{i_{1}}^{i_{s}}g_{1}^{s})(a_{1}^{n})&=
f[(I_{1}^{n}\op{i_{1}}^{i_{s}}g_{1}^{s})\ldots
(I_{n}^{n}\op{i_{1}}^{i_{s}}g_{1}^{s})](a_{1}^{n})
   \\[5pt]
&=f( (I_{1}^{n}\op{i_{1}}^{i_{s}}g_{1}^{s})(a_{1}^{n}),\ldots
,(I_{n}^{n}\op{i_{1}}^{i_{s}}g_{1}^{s})(a_{1}^{n}))    \\[5pt]
&=f( \mu_{1}^{\prime}(\op{i_{1}}^{i_{s}}g_{1}^{s})(a_{1}^{n}),
\ldots ,\mu_{n}^{\prime}(\op{i_{1}}^{i_{s}}g_{1}^{s})(a_{1}^{n})),
\end{array}\vspace{-2mm}
\]
where $\mu_{i}^{\prime}(\op{i_{1}}^{i_{s}}g_{1}^{s})(a_{1}^{n})$
is equal to $a_{i}$ for $i\not\in \{i_{1},\ldots ,i_{s}\}$, and
$\mu_{i}(\op{i_{1}}^{i_{s}}g_{1}^{s})(a_{1}^{n})$ for
$i\in\{i_{1},\ldots ,i_{s}\}$. As $\{i_{1},\ldots
,i_{s}\}=\{j_{1},\ldots ,j_{k}\}$ and the premise of (\ref{13}) is
satisfied, then
\[
\mu_{i}^{\prime }(\op{i_{1}}^{i_{s}}g_{1}^{s})(a_{1}^{n})=\mu
_{i}^{\prime }(\op{j_{1}}^{j_{k}}h_{1}^{k})(a_{1}^{n}),
\]
for all $i=1,\ldots,n$. Hence
\[\arraycolsep=.5mm
\begin{array}{lll}
f\op{i_{1}}^{i_{s}}g_{1}^{s}(a_{1}^{n})&=f(\mu_{1}^{\prime
}(\op{i_{1}}^{i_{s}}g_{1}^{s})(a_{1}^{n}),\ldots ,\mu _{n}^{\prime
}(\op{i_{1}}^{i_{s}}g_{1}^{s})(a_{1}^{n}))   \\[4pt]
&=f(\mu
_{1}^{\prime}(\op{j_{1}}^{j_{k}}h_{1}^{k})(a_{1}^{n}),\ldots
,\mu_{n}^{\prime }(\op{j_{1}}^{j_{k}}h_{1}^{k})(a_{1}^{n}))
=(f\op{j_{1}}^{j_{k}}h_{1}^{k})(a_{1}^{n}),
\end{array}\vspace{-2mm}
\]
which proves the inclusion\footnote{Remind that $f\subset g$ if
and only if pr$_1\,f\subset$ pr$_1\,g$ and $f(x)=g(x)$ for $x\in$
pr$_1\,f$.} $\,f\op{i_{1}}^{i_{s}}g_{1}^{s}\subset
f\op{j_{1}}^{j_{k}}h_{1}^{k}$. The converse inclusion can be
proved analogously. So, the implication (\ref{13}) is valid.
\end{proof}

Basing on the above two propositions we can prove the following
theorem, which was early proved in \cite{Sok} for $n$-ary
operations.

\begin{Theorem}\label{T1}
A $(2,n)$-semigroup $(G;\op{1\;},\ldots,\op{n})$ has a faithful
representation by partial $n$-place functions if and only if for
all $\,g,x_1,\ldots,x_s,y_1,\ldots,y_k\in G$ and
$\,i_1,\ldots,i_s,j_1,\ldots,j_k\in\{1,\ldots,n\}\,$ the following
implication
\begin{equation}\label{14}
\bigwedge\limits_{i=1}^{n}\Big(\mu_i(\opp{i_1}{i_s}x_1^s)
=\mu_i(\opp{j_1}{j_k}y_1^k)\Big)\longrightarrow
g\opp{i_1}{i_s}x_1^s=g\opp{j_1}{j_k}y_1^k
\end{equation}
is satisfied.
\end{Theorem}
\begin{proof}
The necessity of the condition (\ref{14}) follows from Proposition
\ref{P2}. The proof of the sufficiency is based on the
modification on the construction used in the proof of Theorem 3
from \cite{Sok}, where the analogous result was proved for $n$-ary
operations. Let $(G;\op{1\;},\ldots,\op{n})$ be a
$(2,n)$-semigroup and let $G^{\,*}=G\cup\{e_1,\ldots,e_n\}$, where
elements $e_1,\ldots,e_n\not\in G$ are fixed. For
$x_1,\ldots,x_s\in G,$ $\,i_1,\ldots,i_s\in\{1,\ldots,n\}$ and
$\,i=1,\ldots,n,\,$ by $\,\mu_{i}^{*}(\opp{i_1}{i_s}x_1^s)$ we
denote the element of $G^{\,*}$ such that
\[
\mu_{i}^{*}(\opp{i_1}{i_s}x_1^s) =\left\{\begin{array}{cl}
 \mu_{i}(\opp{i_1}{i_s}x_1^s) &\mbox{if}\quad
 i\in\{i_1,\ldots,i_s\},
\\[4pt]
 e_i &\mbox{if}\quad i\not\in\{i_1,\ldots,i_s\}.
\end{array}\right.
\]

Consider the set $A_n= B_{\,n}\cup\{(e_1,\ldots,e_n)\}$, where $B_{\,n}$
is the collection of all $n$-tuples $(x_1,\ldots,x_n)\in
(G^{\,*})^n$ for which there exist $y_1,\ldots,y_s\in G$ and
$i_1,\ldots,i_s\in\{1,\ldots,n\}$ such that
$x_i=\mu_i^*(\opp{i_1}{i_s}y_1^s)$. Next, for every $g\in G$ we
define $n$-place function $\lambda_{g}^{*}:(G^{\,*})^n\rightarrow
G^{\,*}$ putting
\[
\lambda_g^*(x_1^n)=\left\{\begin{array}{cl} g\opp{i_1}{i_s}y_1^s,
& \mbox{if}\quad (x_1,\ldots,x_n)\in B_{\,n},
  \\[4pt]
g & \mbox{if}\quad (x_1,\ldots,x_n)=(e_1,\ldots,e_n).
  \end{array}\right.
\]
For other elements of $(G^{\,*})^n$ this function is not defined.

The mapping $P: g\mapsto \lambda_g^*$ is a faithful representation
of $(G;\op{1\;},\ldots,\op{n})$ by such defined $n$-place
functions. Indeed, if $\,(x_{1},\ldots,x_{n})\in B_{\,n}$, then
\[
\lambda_{g_{1}\op{i\;}g_{2}}^{\ast }(x_{1}^{n})=(g_{1}\op{i\;
}g_{2})\op{i_{1}}^{i_{s}}y_{1}^{s}.
\]
But for $\,i\in \{i,i_{1},\ldots ,i_{s}\}$, we have
\[
\mu _{i}^{\ast }(\op{i\;}g_{2}\op{i_{1}}^{i_{s}}y_{1}^{s})=\mu
_{i}(\op{i\;}g_{2}\op{i_{1}}^{i_{s}}y_{1}^{s})=g_{2}\op{i_{1}}^{i_{s}}y_{1}^{s}=\lambda
_{g_{2}}^{\ast }(x_{1}^{n})
\]
and $\,\mu _{k}^{\ast}
(\op{i\;}g_{2}\op{i_{1}}^{i_{s}}y_{1}^{s})=\mu
_{k}(\op{i_{1}}^{i_{s}}y_{1}^{s})=x_{k} \,$ for all $\,k\neq i.\,$
This means that $\,(x_{1}^{i-1},\lambda_{g_{2}}^{\ast
}(x_{1}^{n}),x_{i+1}^{n})\in B_{\,n}\,$ and in the consequence
\[
g_{1}\op{i\;}g_{2}\op{i_{1}}^{i_{s}}y_{1}^{s}=\lambda
_{g_{1}}^{\ast }(x_{1}^{i-1},\lambda _{g_{2}}^{\ast
}(x_{1}^{n}),x_{i+1}^{n})=\lambda _{g_{1}}^{\ast }\op{i\;}\lambda
_{g_{2}}^{\ast }(x_{1}^{n}).
\]
Thus $\,\lambda_{g_{1}\op{i\;}g_{2}}^{\ast }(x_{1}^{n})=
\lambda_{g_{1}}^{\ast }\op{i\;}\lambda_{g_{2}}^{\ast
}(x_{1}^{n})\,$ for all $\,(x_1,\ldots,x_n)\in B_{\,n}$.

In the case $(x_{1},\ldots,x_{n})=(e_{1},\ldots ,e_{n})$ we have
$\lambda _{g_{1}\op{i\;}g_{2}}^{\ast
}(e_{1}^{n})=g_{1}\op{i\;}g_{2}$ and $\lambda_{g_{2}}^{\ast
}(e_{1}^{n})=g_{2}=\mu _{i}(\op{i\;}g_{2})= \mu _{i}^{\ast
}(\op{i\;}g_{2})$. Since $\mu_{k}^{\ast }(\op{i\;}g_{2})=e_{k}$
for $k\neq i$, the element $(e_{1}^{i-1},\lambda _{g_{2}}^{\ast
}(e_{1}^{n}),e_{i+1}^{n})$ lies in $B_{\,n}$ and $\lambda
_{g_{1}}^{\ast }(e_{1}^{i-1},\lambda _{g_{2}}^{\ast
}(e_{1}^{n}),e_{i+1}^{n})=g_{1}\op{i\;}g_{2}$. So,
$\,\lambda_{g_{1}\op{i\;}g_{2}}^{\ast }(e_{1}^{n})=
\lambda_{g_{1}}^{\ast }\op{i\;}\lambda_{g_{2}}^{\ast
}(e_{1}^{n})$. Thus, in any case we have
$P(g_{1}\op{i\;}g_{2})=P(g_{1})\op{i\;}P(g_{2})$, which completes
the proof that $P$ is a homomorphism. Hence $P$ is a
representation of $(G;\op{1\;},\ldots,\op{n})$ by partial
$n$-place functions.

This representation is faithful because $P(g_1)=P(g_2)$, i.e.
$\lambda_{g_{1}}^{\ast }=\lambda _{g_{2}}^{\ast }$, implies
$\lambda _{g_{1}}^{\ast }(e_{1}^{n})=\lambda _{g_{2}}^{\ast
}(e_{1}^{n})$, whence $g_{1}=g_{2}$.
\end{proof}

\begin{Theorem}\label{T2}
Every $(2,n)$-semigroup satisfying the condition $(\ref{14})$ has
a faithful representation by full $n$-place functions.
\end{Theorem}
\begin{proof} Let $(G;\op{1\;},\ldots,\op{n})$ be some
$(2,n)$-semigroup. By Theorem~\ref{T1} it is isomorphic to some
$(2,n)$-semigroup $(\Phi;\op{1\;},\ldots,\op{n})$ of partial
$n$-place functions, where $\Phi\subset\mathcal{F}(A^n,A)$.
Consider now the set $A_0=A\cup\{c\}$, where $c\not\in A$, and the
extension $f^{0}$ of $f\in\Phi$ defined in the following way:
$$
f^{0}(x_1^n)=\left\{\begin{array}{cl}
 f(x_0^n)& \mbox{if}\quad (x_1,\ldots,x_n)\in\pr f, \\[4pt]
  c & \mbox{if}\quad (x_1,\ldots,x_n)\not\in\pr f
\end{array}\right.$$
for all $x_1,\ldots,x_n\in A_0$. It is clear that $f^{0}$ is a
full $n$-place function on $A_0$.

Let us show that the mapping $f\mapsto f^{0}$ is an isomorphism of
$(\Phi;\op{1\;},\ldots,\op{n})$ onto
$(\Phi_0;\op{1\;},\ldots,\op{n})$, where
$\Phi_0=\{f^{0}\,|\,f\in\Phi\}$. Indeed, if $f,g\in\Phi$, then in
the case $(x_1,\ldots,x_n)\in\pr(f\op{i\;}g)$ we have
$(f\op{i\;}g)^0(x_1^n)=f\op{i\;}g(x_1^n)$ and
$(x_1,\ldots,x_n)\in\pr g$,
$\,(x_1^{i-1},g(x_1^n),x_{i+1}^n)\in\pr f.$ Hence
$g^{0}(x_1^n)=g(x_1^n)$, and in the consequence
\[\arraycolsep=.5mm
\begin{array}{rl}
f^{0}\op{i\;}g^{0}(x_1^n)&=f^{0}(x_1^{i-1},g^{0}(x_1^n),x_{i+1}^n)=
f^{0}(x_1^{i-1},g(x_1^n),x_{i+1}^n) \\[4pt]
&=f(x_1^{i-1},g(x_1^n),x_{i+1}^n)=f\op{i\;}g(x_1^n).
\end{array}
\]
Thus $(f\op{i}g)^{0}(x_1^n)=f^{0}\op{i\;}g^{0}(x_1^n)$.

In the case when $(x_1,\ldots,x_n)\not\in\pr(f\op{i\;}g)$ we have
$(f\op{i\;}g)^0(x_1^n)=c$ and $(x_1,\ldots,x_n)\not\in\pr g$, or
$(x_1,\ldots,x_n)\in\pr g$ and
$(x_1^{i-1},g(x_1^n),x_{i+1}^n)\not\in\pr f$. If
$(x_1,\ldots,x_n)\not\in\pr g$, then $g^{0}(x_1^n)=c$ and
$$
f^{0}\op{i\;}g^{0}(x_1^n)=f^{0}(x_1^{i-1},g^{0}(x_1^n),x_{i+1}^n)=
f^{0}(x_1^{i-1},c,x_{i+1}^n)=c.
$$
If $(x_1,\ldots,x_n)\in\pr g$ and $(x_1^{i-1},g(x_1^n),x_{i+1}^n)\not\in\pr f$,
then $g^{0}(x_1^n)=g(x_0^n)$ and
$$
f^{0}\op{i\;}g^{0}(x_1^n)=f^{0}(x_1^{i-1},g^{0}(x_1^n),x_{i+1}^n)
=f^{0}(x_1^{i-1},g(x_1^n),x_{i+1}^n)=c.
$$
So, in all cases $(f\op{i\;}g)^0(x_1^n)=f^{0}\op{i\;}g^{0}(x_1^n)$
for all $(x_1,\ldots,x_n)\in A_0^n$. Thus, the mapping $f\mapsto
f^{0}$ is a homomorphism.

If $f^{0}=g^{0}$ for some $f,g\in\Phi$, then $\pr f=\pr g$ and
$f(x_1^n)=g(x_1^n)$ for all $x_1^n\in\pr f,$ i.e. $f=g$. So,
$f\mapsto f^{0}$ is an isomorphism.
\end{proof}

\begin{Definition}\rm
A $(2,n)$-semigroup $(G;\op{1\;},\ldots,\op{n})$ is called {\it
unitary}, if it contains {\it selectors}, i.e. such elements
$e_1,\ldots, e_n$ for which the equalities
\begin{eqnarray}\label{15}
&&g\op{i}e_i=e_i\op{i}g=g,\\[4pt]
 \label{16}
&&e_k\op{i}g=e_k
\end{eqnarray}
are valid for all $g\in G$ and $i,k\in\{1,\ldots,n\}$, where $i\ne
k$.
\end{Definition}

Let $(G;\op{1\;},\ldots,\op{n})$ be a $(2,n)$-semigroup. We say
that a unitary $(2,n)$-semigroup
$(G^{\,*},\op{1\;},\ldots,\op{n})$ with selectors $e_1,\ldots,e_n$
is a {\it unitary extension} of a $(2,n)$-semigroup
$(G;\op{1\;},\ldots,\op{n})$, if
\begin{itemize}
  \item[(a)] $G\subset G^{\,*}$,
  \item[(b)] $G\cap\{e_1,\ldots,e_n\}=\varnothing$,
  \item[(c)] $G\cup\{e_1,\ldots,e_n\}$ is a generating set of a $(2,n)$-semigroup
$(G^{\,*},\op{1\;},\ldots,\op{n})$.
\end{itemize}
%with adjoining complete collection of selectors.

\begin{Theorem}\label{T3}
Every representable $(2,n)$-semigroup can be isomorphically
embedded into a unitary extension of some $(2,n)$-semigroup.
\end{Theorem}
\begin{proof} By Theorem~\ref{T1} every representable $(2,n)$-semigroup
$(G;\op{1\;},\ldots,\op{n})$ satisfies the condition (\ref{14}).
By Theorem~\ref{T2} it is isomorphic to some $(2,n)$-semigroup
$(\Phi_0;\op{1\;},\ldots,\op{n})$ of full $n$-place functions on
some set $A_0$.

Consider the set $\{I^n_1,\ldots,I^n_n\}$ of $n$-place projectors
on $A_0$ and the family of subsets $\,(F_k(\Phi_0))_{k\in\mathbb
N}\,$ satisfying the following two conditions:
\begin{enumerate}
  \item[1)] $F_0(\Phi_0)=\Phi_0\cup\{I^n_1,\ldots,I^n_n\}$,
  \item[2)]  $f,g\in F_k(\Phi_0)\longrightarrow f\op{i\;}g\in
  F_{k+1}(\Phi_0)$.
\end{enumerate}
Obviously $\{I^n_1,\ldots,I^n_n\}\subset F_k(\Phi_0)\,$ for every
$k\in\mathbb N$. Moreover, if $f\in F_k(\Phi_0)$, then
$f=f\op{i\;}I^n_i\in F_{k+1}(\Phi_0)$. So, $\Phi_0\subset
F_k(\Phi_0)\subset F_{k+1}(\Phi_0)$ for every $k\in\mathbb N.$

Now let $\Phi^{*}=\bigcup\limits_{k=0}^{\infty}F_k(\Phi_0)$. Then
$\Phi_0\cup\{I^n_1,\ldots,I^n_n\}=F_0(\Phi_0)\subset\Phi^{*}$ and
$\{I^n_1,\ldots,I^n_n\}\cap\Phi_0=\varnothing$. If $f\in
F_n(\Phi_0)$, $g\in F_m(\Phi_0)$ for some $n,m\in\mathbb{N}$, then
$f,g\in F_k(\Phi_0)$, where $k=\max(n,m)$. Therefore
$f\op{i\;}g\in F_{k+1}(\Phi_0)\subset\Phi^{*}$. This means that
the set $\Phi^{*}$ is closed with respect to the operations
$\op{1\;},\ldots,\op{n}$ and contains the selectors
$I^n_1,\ldots,I^n_n$. Hence $(\Phi^{*};\op{1\;},\ldots,\op{n})$ is
a unitary extension of a $(2,n)$-semigroup
$(\Phi_0;\op{1\;},\ldots,\op{n})$.
\end{proof}

Let $(P_i)_{i\in I}$ be the family of representations of a
$(2,n)$-semigroup $(G;\op{1\;},\ldots,\op{n})$ by $n$-place
functions defined on sets $(A_i)_{i\in I}$, respectively. By the
{\em union} of this family we mean the mapping $P:g\mapsto P(g)$,
where $g\in G$, and $P(g)$ is an $n$-place function on
$A=\bigcup\limits_{i\in I}A_i$ defined by
\[
P(g)=\bigcup\limits_{i\in I}P_i(g)\, .
\]
If $A_i\cap A_j=\emptyset$ for all $i,j\in I$, $i\neq j$, then $P$
is called the {\em sum} of $(P_i)_{i\in I}$ and is denoted by
$P=\sum_{i\in I}P_i$. It is not difficult to see that the sum of
representations is a representation, but the union of
representations may not be a representation.

Let $(G;\op{1\;},\ldots,\op{n})$ be a $(2,n)$-semigroup. A binary
relation $\rho\subset G\times G$ is called
\begin{itemize}
\item \textit{$v$-regular}, if
\[
\bigwedge\limits_{i=1}^{n}(x_i,y_i)\in\rho\longrightarrow
(g\opp{i_1}{i_s}u_1^s\, ,\,g\opp{j_1}{j_k}v_1^k)\in\rho
\]
for all $g\in G$ and $x_i=\mu_i(\opp{i_1}{i_s}u_1^s)$,
$\,y_i=\mu_i(\opp{j_1}{j_k}v_1^k)$, $i=1,\ldots,n$, where
$u_1,\ldots,u_s\in G,$ $\,v_1,\ldots,v_k\in G$,
\item \textit{$l$-regular}, if
\[
(x,y)\in\rho\longrightarrow(x\op{i\;}z\, ,\,y\op{i\;}z)\in\rho
\]
for all $x,y,z\in G,$ $\,i=1,\ldots,n$,
\item \textit{$v$-negative}, if
\[
\Big(x\opp{i_1}{i_s}y_1^n\, ,\,
\mu_i(\opp{i_1}{i_s}y_1^s)\Big)\in\rho
\]
for all $i\in\{i_1,\ldots,i_s\}$ and $x,y_1,\ldots,y_s\in G.$
\end{itemize}

A nonempty subset $W$ of $G$ is called an \textit{$l$-ideal}, if
the implication
\[
g\opp{i_1}{i_s}x_1^s\not\in
W\longrightarrow\mu_i(\opp{i_1}{i_s}x_1^s)\not\in W
\]
is valid for all $g,x_1,\ldots,x_s\in G$ and
$i\in\{i_1,\ldots,i_s\}\subset\{1,\ldots,n\}$.

\bigskip

By a \textit{determining pair} of a $(2,n)$-semigroup
$(G;\op{1\;},\ldots,\op{n})$ we mean an ordered pair
$(\mathcal{E},W)$, where $\mathcal{E}$ is a symmetric and
transitive binary relation defined on a unitary extension
$(G^{\,*};\op{1\;},\ldots,\op{n})$ of a $(2,n)$-semigroup
$(G;\op{1\;},\ldots,\op{n})$ and $W$ is a subset of $G^{\,*}$ such
that
\begin{enumerate}
\item[1)] $G\cup\{e_1,\ldots,e_n\}\subset\pr\mathcal{E}$,
\item[2)] $\{e_1,\ldots,e_n\}\cap W=\varnothing$,
\item[3)] $\bigwedge\limits_{i=1}^{n}e_i\equiv
x_i(\mathcal{E})\longrightarrow g\equiv
g\opp{i_1}{i_s}y_1^s(\mathcal{E})$, where
$x_i=\mu_i(\opp{i_1}{i_s}y_1^s)$, $i=1,\ldots,n$,
$\,y_1,\ldots,y_s\in G^{\,*}$,
\item[4)] $\bigwedge\limits_{i=1}^{n}x_i\equiv y_i(\mathcal{E})\longrightarrow
g\opp{i_1}{i_s}u_1^s\equiv g\opp{j_1}{j_k}v_1^k(\mathcal{E})\,$
for all $g\in G,$  $x_i=\mu_i(\opp{i_1}{i_s}u_1^s)\in G,\,$
$y_i=\mu_i(\opp{j_1}{j_k}v_1^k)$, $\,i=1,\ldots,n$,  where
$u_1,\ldots,u_s\in G^*,$ $v_1,\ldots,v_k\in G^*,$
\item[5)] if $W\neq\varnothing$, then $W$ is an $\mathcal{E}$-class and
$W\cap G$ is an $l$-ideal of $G$.
\end{enumerate}

Let $\left(H_a\right)_{a\in A}$ be a collection of
$\mathcal{E}$-classes (uniquely indexed by elements of $A$) such
that $H_a\neq W$ and
$H_a\cap(G\cup\{e_1,\ldots,e_n\})\neq\emptyset$ for all $a\in A$.
Consider the set $\frak{A}$ of elements $(a_1,\ldots,a_n)\in A^n$
satisfying one from the following conditions
\begin{itemize}
\item[$(a)$]
$\,H_{a_i}=\mathcal{E}\langle\mu^*(\opp{i_1}{i_s}y_1^s)\rangle\,$
for all $\,i=1,\ldots,n\,$ and some $\,y_1,\ldots,y_s\in
G$,\footnote{\,$\mathcal{E}\langle x\rangle$ denotes the
$\mathcal{E}$-class of $x$.}
\item[$(b)$] $\,H_{a_i}=\mathcal{E}\langle e_i\rangle\,$ for all $\,i=1,\ldots,n$.
\end{itemize}

Next, using this set, for each element $g\in G$, we define an
$n$-place function $P_{(\mathcal{E},W)}(g)$ on $A$ putting
\[
(a_1^n,b)\in P_{(\mathcal{E},W)}(g)\longleftrightarrow
(a_1,\ldots,a_n)\in\frak{A}\,\wedge\left\{\begin{array}{rll}
\!\!g\opp{i_2}{i_s}y_1^s\in H_b&{\rm if\;holds }\;\;(a),\\[6pt]
\!g\in H_b&{\rm if\;holds}\;\;(b).
\end{array}\right.
\]
\smallskip
\begin{Proposition}\label{P3}
If a $(2,n)$-semigroup $(G;\op{1\;},\ldots,\op{n})$ is
representative, then a mapping $\,g\mapsto
P_{(\mathcal{E},W)}(g)$, where $(\mathcal{E},W)$ is its
determining pair, is a representation of this $(2,n)$-semigroup by
$n$-place functions.
\end{Proposition}
\begin{proof} We must show that
\begin{equation}\label{17}
P_{(\mathcal{E},W)}(g_1\op{i\;}g_2)=P_{(\mathcal{E},W)}(g_1)\op{i\;}
P_{(\mathcal{E},W)}(g_2)
\end{equation}
for all $g_1,g_2\in G$ and $i=1,\ldots,n$.

Let $(a_1^n,b)\in P_{(\mathcal{E},W)}(g_1\op{i\;}g_2)$ for some
$i=1,\ldots,n$. Then, according to the definition of
$P_{(\mathcal{E},W)}$, we have $(a_1,\ldots,a_n)\in\frak{A}$. In
the case $(a)$, we have also
$(g_1\op{i\;}g_2)\opp{i_1}{i_s}y_1^s\in H_b$. But $H_b\cap
W=\varnothing$, so, $(g_1\op{i\;}g_2)\opp{i_1}{i_s}y_1^s\not\in W.$
Therefore $\mu_i^*(\op{i\;}g_2\opp{i_1}{i_s}y_1^s)\not\in W$, i.e.
$g_2\opp{i_1}{i_s}y_1^s\not\in W$, because $W$ is an $l$-ideal.
Assume that $g_2\opp{i_1}{i_s}y_1^s\in H_c$. Since\vspace{-4mm}
\[
\mu_k^*(\op{i}g_2\opp{i_1}{i_s}y_1^s)=\mu_k^*(\opp{i_1}{i_s}y_1^s)\not\in
W\quad\mbox{for}\quad k\neq i,
\]
we have $(a_1^{i-1},c,a_{i+1}^n)\in\frak{A}$, which, together with
$(a_1^n,c)\in P_{(\mathcal{E},W)}(g_2)$, implies $(a_1^n,b)\in
P_{(\mathcal{E},W)}(g_1)\op{i\;}P_{(\mathcal{E},W)}(g_2)$.

In the case $(b)$ we obtain $g_1\op{i\;}g_2\in H_b$. Thus
$g_2\not\in W,$ and in the consequence, $g_2\in H_c$ for some
$c\in A$. Therefore $(a_1^{i-1},c,a_{i+1}^n)\in\frak{A}$, whence
\[
(a_1^{i-1}\,c\,a_{i+1}^n,b)\in
P_{(\mathcal{E},W)}(g_1)\qquad\mbox{and}\qquad(a_1^n,c)\in
P_{(\mathcal{E},W)}(g_2).
\]
Consequently $\,(a_1^n,b)\in
P_{(\mathcal{E},W)}(g_1)\op{i\;}P_{(\mathcal{E},W)}(g_2)$, which
shows that the inclusion
\[
P_{(\mathcal{E},W)}(g_1\op{i\;}g_2)\subset
P_{(\mathcal{E},W)}(g_1)\op{i\;}P_{(\mathcal{E},W)}(g_2)
\]
is valid in any case.

Now let $(a_1^n,b)\in
P_{(\mathcal{E},W)}(g_1)\op{i\;}P_{(\mathcal{E},W)}(g_2)$. Then
there exists $c\in A$ such that $\,(a_1^n,c)\in
P_{(\mathcal{E},W)}(g_2)\,$ and $\,(a_1^{i-1}c\,a_{i+1}^n,b)\in
P_{(\mathcal{E},W)}(g_1)$.

If $H_{a_i}=\mathcal{E}\langle e_i\rangle$ for all $i=1,\ldots,n$,
then $g_2\in H_c$. Thus $(a_1^{i-1},c,a_{i+1}^n)\in\frak{A}$ and
$g_1\op{i\;}g_2\in H_b$, whence $(a_1^n,b)\in
P_{(\mathcal{E},W)}(g_1\op{i\;}g_2)$.

If $H_{a_i}=\mathcal{E}\langle\mu^*(\opp{i_1}{i_s}y_1^s)\rangle$
for all $i=1,\ldots,n$ and some $y_1,\ldots,y_s\in G,$ then
$\,H_b=\mathcal{E}\langle
g_1\op{i\;}g_2\opp{i_1}{i_s}y_1^s\rangle$, because
$\,H_c=\mathcal{E}\langle g_2\opp{i_1}{i_2}y_1^s\rangle=
\mathcal{E}\langle\mu_i^*(\op{i\;}g_2\opp{i_1}{i_s}y_1^s)\rangle$
and \vspace{-2mm}
\[
H_k=\mathcal{E}\langle\mu_i^*(\opp{i_1}{i_s}y_1^s)\rangle=
\mathcal{E}\langle\mu_i^*(\op{i\;}g_2\opp{i_1}{i_s}y_1^s)\rangle,\quad
k\neq i.
\]
Therefore $\,(a_1^n,b)\in P_{(\mathcal{E},W)}(g_1\op{i\;}g_2)$.
So,
\[
P_{(\mathcal{E},W)}(g_1)\op{i\;}P_{(\mathcal{E},W)}(g_2)\subset
P_{(\mathcal{E},W)}(g_1\op{i\;}g_2),
\]
which, together with the previous inclusion, proves (\ref{17}).
\end{proof}

According to \cite{ST}, the representation $P_{(\mathcal{E},W)}$,
uniquely determined by the pair $(\mathcal{E},W)$, will be called
\textit{simplest}.

\begin{Theorem}\label{T4}
Any representation of a $(2,n)$-semigroup by $n$-place functions
is a union of some family of its simplest representations.
\end{Theorem}
\begin{proof} Let $P$ be a representation of a $(2,n)$-semigroup
$(G;\op{1\;},\ldots,\op{n})$ by $n$-place functions defined on
$A$, and let $c\not\in A$ be some fixed element. For every $g\in
G$ we define on $A_0=A\cup\{c\}$ an $n$-place function $P^*(g)$
putting:
\[
P^*(g)=\left\{\begin{array}{cl}
  P(g)(a_1^n) & \mbox{if}\quad (a_1,\ldots,a_n)\in\pr P(g), \\[4pt]
  c & \mbox{if}\quad (a_1,\ldots,a_n)\not\in\pr P(g).
\end{array}\right.
\]
It is not difficult to see that $P^*$ is a representation of
$(G;\op{1\;},\ldots,\op{n})$ by $n$-place functions defined on
$A_0$, and $P(g)\mapsto P^*(g)$, where $g\in G$, is an isomorphism
of $(P(G);\op{1\;},\ldots,\op{n})$ onto
$(P^*(G);\op{1\;},\ldots,\op{n})$. Because
$G\cup\{e_1,\ldots,e_n\}$ is a generating set of a unitary
extension $(G^*;\op{1\;},\ldots,\op{n})$ with selectors
$e_1,\ldots,e_n$, then putting $P^*(e_i)=I_i^n$, $i=1,\ldots,n$,
where $I_i^n$ is the $i$-th $n$-place projector of $A_0$, we
obtain a unique extension of $P^*$ from $G$ to $G^*$.

For any $(a_1,\ldots,a_n)\in A^n$ we define on $G^*$ an
equivalence $\Theta_{a_1^n}$ such that
\[
x\equiv y(\theta_{a_1^n})\longleftrightarrow
P^*(x)(a_1^n)=P^*(y)(a_1^n).
\]
This equivalence is $v$-regular. Indeed, if $x_i\equiv
y_i(\theta_{a_1^n})$, i.e. $P^*(x_i)(a_1^n)=P^*(y_i)(a_1^n)$ for
$i=1,\ldots,n$, then for $g\in G$ and
$x_i=\mu_i(\opp{i_1}{i_s}u_1^s)$,
$y_i=\mu_i(\opp{j_1}{j_k}v_1^k)$, $i=1,\ldots,n$, by
Proposition~\ref{P1}, we have
\[\arraycolsep=.5mm
\begin{array}{rl}
P^*(g\opp{v_1}{i_s}u_1^s)(a_1^n)&
=P^*(g)\op{i_1}P^*(u_1)\op{i_2}\ldots\op{i_s}P^*(u_s)(a_1^n) \\[8pt]
&=P^*(g)\Big[\Big(I^n_1\op{i_1}P^*(u_1)\op{i_2}\ldots\op{i_s}P^*(u_s)\Big)\ldots \\[8pt]
&\hspace*{25mm}\ldots
\Big(I^n_n\op{i_1}P^*(u_1)\op{i_2}\ldots\op{i_s}P^*(u_s)\Big)\Big](a_1^n)
 \\[6pt]
&=P^*(g)\Big(P^*(\mu_1(\opp{i_1}{i_s}u_1^s))(a_1^n)\, ,\;\ldots\;
, \,P^*(\mu_n(\opp{i_1}{i_s}u_1^s))(a_1^n)\Big) \\[6pt]
&=P^*(g)\Big(P^*(x_1)(a_1^n)\, ,\;\ldots\; ,\,P^*(x_n)(a_1^n)\Big)
\\[6pt]
&=P^*(g)\Big(P^*(y_1)(a_1^n)\, ,\;\ldots\; ,\,
P^*(y_n)(a_1^n)\Big)\\[6pt]
&=\ldots=P^*(g\opp{j_1}{j_k}v_1^k)(a_1^n). \end{array}
\]
So, $g\opp{i_1}{i_s}u_1^s\equiv
g\opp{j_1}{j_k}v_1^k(\theta_{a_1^n})$. This proves the
$v$-regularity of the equivalence $\theta_{a_1^n}$. All subsets of
the form
\[
H_{b}^{a_1^n}=\{x\in G^*\,|\,(a_1^n,b)\in P^*(x)\}
\]
are, of course, the equivalence classes of this relation.
Moreover, the pair $(\mathcal{E}_{a_1^n},W_{a_1^n})$, where
\begin{eqnarray*}
&\mathcal{E}_{a_1^n}=\theta_{a_1^n}\cap\theta_{a_1^n}
\Big(G\cup\{e_1,\ldots,e_n\}\Big)\times\theta_{a_1^n}\Big(G\cup\{e_1,\ldots,e_n\}
\Big), \\[4pt]
&W_{a_1^n}=\{x\in G^*\,|\,P^*(x)(a_1^n)=c \},
\end{eqnarray*}
is a determining pair of a $(2,n)$-semigroup
$(G;\op{1\;},\ldots,\op{n})$.

We prove that a representation $P$ is a union of the family of
simplest representations $P_{a_1^n}$ of
$(G;\op{1\;},\ldots,\op{n})$ induced by a determining pair
$(\mathcal{E}_{a_1^n},W_{a_1^n})$, i.e. that
\begin{equation}\label{18}
P(g)=\bigcup\limits_{a_1^n\in A^n}\!\!\! P_{a_1^n}(g)
\end{equation}
for every $\,g\in G$. Indeed, if $(b_1^n,d)\in P(g)$, then $g\in
H_{d}^{b_1^n}$, where $b_1^n\in A^n,$ $\,d\in A$. But $e_i\in
H_{b_i}^{b_1^n}$, $\,i=1,\ldots,n$, and $g\in H_{d}^{b_1^n}$,
imply, according to the definition, $(b_1^n,d)\in P_{b_1^n}(g)$.
Therefore $\,(b_1^n,d)\in\bigcup\limits_{a_1^n\in
A^n}P_{a_1^n}(g)$. \ So,
\[
  P(g)\subset\bigcup\limits_{a_1^n\in A^n}\!\!\! P_{a_1^n}(g).
\]

Conversely, if $(b_1^n,d)\in P_{a_1^n}(g)$ for some $a_1^n\in
A^n$, then $g\opp{i_1}{i_s}y_1^s\in H_{d}^{a_1^n}$ (for
$H_{b_i}^{a_1^n}=\mathcal{E}_{a_1^n}\langle\mu_i^*(\opp{i_1}{i_s}y_1^s)\rangle$,
$\,i=1,\ldots,n$), or $\,g\in H_d^{a_1^n}$ (for
$H_{b_i}^{a_1^n}=\mathcal{E}_{a_1^n}\langle e_i\rangle$,
$\,i=1,\ldots,n$, where $b_1^n\in A^n$).

For $\mu_i^*(\opp{i_1}{i_s}y_1^s)\in H_{b_i}^{a_1^n}$ we have
$b_i=P\Big(\mu_i^*(\opp{i_1}{i_s}y_1^s)\Big)(a_1^n)$,
$\,i=1,\ldots,n$. From $g\opp{i_1}{i_s}y_1^s\in H_{d}^{a_1^n}$ we
obtain $d=P\Big(g\opp{i_1}{i_s}y_1^s\Big)(a_1^n)$. But $P$ is a
homomorphism, so,
\[
d=P(g)\left(P\Big(\mu_1^*(\opp{i_1}{i_s}y_1^s)\Big)(a_1^n),\ldots,
P\Big(\mu_n^*(\opp{i_1}{i_s}y_1^s)\Big)(a_1^n)\right)=
P(g)(b_1^n).
\]
Hence $\,(b_1^n,c)\in P(g)$.

For $e_i\in H_{b_i}^{a_1^n}$, $\,i=1,\ldots,n$, we get
$(a_1^n,b_i)\in P^*(e_i)=I^n_i$, $i=1,\ldots,n$, whence $a_i=b_i$
for all $\,i=1,\ldots,n$. So, $\,g\in H_{d}^{b_1^n},$ \ i.e.
$(b_1^n,d)\in P(g)$.

Thus, in these both cases we have
\[
\bigcup\limits_{a_1^n\in A^n}\!\!\!P_{a_1^n}(g)\subset P(g).
\]
which, together with the previous inclusion, proves (\ref{18}).
\end{proof}

\medskip

Let $P$ be a representation of a $(2,n)$-semigroup
$(G;\op{1\;},\ldots,\op{n})$ by $n$-place functions. Define on $G$
two binary relations $\chi_{_P}$ and $\varepsilon_{_P}$ putting
\begin{eqnarray*}
&&(g_1,g_2)\in\chi_{_P}\longleftrightarrow\pr
P(g_1)\subset\pr P(g_2),\\[4pt]
&&(g_1,g_2)\in\varepsilon_{_P}\longleftrightarrow P(g_1)=P(g_2).
\end{eqnarray*}
It is not difficult to see that the relation  $\chi_{_P}$ is
reflexive and transitive, i.e. it is a quasi-order. The relation
$\varepsilon_{_P}$ is an equivalence on $G$. If a representation
$P$ is faithful, then
$\varepsilon_{_P}=\triangle_G=\{(g,g)\,|\,g\in G\}$. In the case
when $P$ is a representation by full $n$-place functions, we have
$\chi_{_P}=G\times G$. Moreover, if $P$ is a sum of a family
$(P_i)_{i\in I}$ of representations $P_i$, then
\begin{equation}\label{19}
 \chi_{_P}=\bigcap\limits_{i\in
I}\chi_{_{P_i}}\quad{\rm and }\quad
\varepsilon_{_P}=\bigcap\limits_{i\in I}\varepsilon_{_{P_i}}.
\end{equation}

Following \cite{ST} an algebraic system
$(\Phi;\op{1\,},\ldots,\op{n},\chi_{_\Phi})$, where
$(\Phi;\op{1\,},\ldots,\op{n})$ is a $(2,n)$-semigroup of
$n$-place functions on a set $A$ and
\[
\chi_{_\Phi}=\{(f,g)\in\Phi\times\Phi\,|\,\pr f\subset\pr g\},
\]
is called a \textit{projection quasi-ordered $(2,n)$-semigroup of
$n$-place functions.} It is characterized by the following
theorem.

\begin{Theorem}\label{T5}
An algebraic system $(G;\op{1\,},\ldots,\op{n},\chi)$, where
$(G;\op{1\,},\ldots,\op{n})$ is a $(2,n)$-semigroup, $\chi$ is a
binary relation on $G$, is isomorphic to a projection
quasi-ordered  $(2,n)$-semigroup of $n$-place functions if an only
if it satisfies the condition $(\ref{14})$ and $\chi$ is an
$l$-regular, $v$-negative quasi-order.
\end{Theorem}
\begin{proof} {\it Necessity}. Let
$(\Phi;\op{1},\ldots,\op{n},\chi_{_\Phi})$ be a projection
quasi-ordered $(2,n)$-semigroup of $n$-place functions. It is
clear, that the relation $\chi_{_\Phi}$ is reflexive and
transitive, i.e. it is a quasi-order. By Theorem~\ref{T1} the
condition (\ref{14}) is satisfied.

Assume that for some $f,g\in\Phi$ we have $(f,g)\in\chi_{_\Phi}$,
i.e. $\pr f\subset\pr g$. If $(a_1,\ldots,a_n)\in\pr(f\op{i\,}h)$,
where $h\in\Phi$, then there exists $c\in A$ such that
$(a_1^n,c)\in f\op{i\,}H$, whence $(a_1^n,b)\in h$ and
$(a_1^{i-1}\,b\,a_{i+1}^n,c)\in f$ for some $b\in A$. So,
$(a_1^{i-1}\,b\,a_{i+1}^n)\in\pr f.$ Therefore
$(a_1^{i-1}\,b\,a_{i+1}^n)\in\pr g$, whence
$(a_1^{i-1}\,b\,a_{i+1}^n,d)\in g$ for some $d\in A.$ Thus,
$(a_1^n,b)\in h$ and $(a_1^{i-1}\,b\,a_{i+1}^n,d)\in g.$ Hence
$(a_1^n,d)\in g\op{i\,}h$, i.e.
$(a_1,\ldots,a_n)\in\pr(g\op{i\,}h)$, which proves the inclu\-sion
$\pr(f\op{i\,}h)\subset\pr(g\op{i\,}h)$. So, the relation
$\chi_{_\Phi}$ is $l$-regular.

Now let $(a_1,\ldots,a_n)\in\pr(f\opp{i_1}{i_s}g_1^s)$ for some
$f,g_1,\ldots,g_s\in\Phi$. This means that
$f\opp{i_1}{i_s}g_1^s\langle
a_1^n\rangle\neq\varnothing$.\footnote{\,The expression $f\langle
a_1^n\rangle$ denotes the set $\{f(a_1^n)\}$, if
$(a_1,\ldots,a_n)\in\mathrm{pr}_{1}f$, and empty set, if
$(a_1,\ldots,a_n)\not\in\mathrm{pr}_{1}f$.} Thus, according to
(\ref{12}), we obtain
\[\arraycolsep=.5mm
\begin{array}{rl}
\varnothing&\neq f\opp{i_1}{i_s}g_1^s\langle a_1^n
\rangle=f[(I^n_1\opp{i_1}{i_s}g_1^s)\ldots
(I^n_n\opp{i_1}{i_s}g_1^s)]\langle a_1^n\rangle\\[3mm]
&=f\Big(I^n_1\opp{i_1}{i_s}g_1^s\langle a_1^n\rangle\,
,\ldots,\,I^n_n\opp{i_1}{i_s}g_1^s\langle a_1^n\rangle\Big),
\end{array}
\]
where $I^n_1,\ldots,I^n_n$ are $n$-place projectors on $A$. Hence
$I^n_i\opp{i_1}{i_s}g_1^s\langle a_1^n\rangle\neq\varnothing$ for
any $i\in\{i_1,\ldots,i_s\}$. Thus
$\mu_i(\opp{i_1}{i_s}g_1^s)\langle a_1^n\rangle\neq\varnothing$,
i.e. $(a_1,\ldots,a_n)\in\pr\mu_i(\opp{i_1}{i_s}g_1^s)$. This
shows that the relation $\chi_{_\Phi}$ is $v$-negative.

{\it Sufficiency}. Let all the conditions of the theorem will be
satisfied by an algebraic system $(G;\op{1\,},\ldots,\op{n},\chi)$
and let $G^*=G\cup\{e_1,\ldots,e_n\}$, where
$e_1,\ldots,e_n\not\in G$. Consider the set $\mathbf{B_0}$ defined
in the following way:
\[
(x_1,\ldots,x_n)\in\mathbf{B_0}\longleftrightarrow(\forall
i=1,\ldots,n)\ x_i=\mu_i^*(\opp{i_1}{i_s}y_1^s)
\]
for some $y_1,\ldots,y_s\in G$, $i_1,\ldots,i_s\in\{1,\ldots,n\}$,
where $\mu_i^*(\opp{i_1}{i_s}y_1^s)$ denotes an element from $G^*$
defined in the proof of Theorem~\ref{T1}.

Let $a\in G$ be fixed. For every $g\in G$ we define an $n$-place
function $P_a(g)$ from
$\mathbf{B}=\mathbf{B_0}\cup\{(e_1,\ldots,e_n)\}$ to $G$ putting
\[
P_a(g)(x_1^n)=\left\{
\begin{array}{cl}
g\opp{i_1}{i_s}y_1^s & \mbox{if }\;
x_i=\mu_i^*(\opp{i_1}{i_s}y_1^s),\; i=1,\ldots,n,\mbox{ and }
a\sqsubset
g\opp{i_1}{i_s}y_1^s\\
&\quad\;\mbox{for some}\;\; y_1,\ldots,y_s\in G,\\[4pt]
   g&\mbox{if }\;(x_1,\ldots,x_n)=(e_1,\ldots,e_n)\mbox{ and } a\sqsubset g,
\end{array}\right.
\]
where $(x_1,\ldots,x_n)\in\mathbf{B}$ and $a\sqsubset g
\longleftrightarrow (a,g)\in\chi$. It is clear that $P_a(g)$ is a
partial $n$-place function on $G^*$.

Let us show that $P_a:g\mapsto P_a(g)$ is a homomorphism, i.e. we
verify the identity
\begin{equation}\label{20}
 P_a(g_1\op{i\,}g_2)=P_a(g_1)\op{i\,}P_a(g_2).
\end{equation}

For this consider $g_1,g_2\in G$, $(x_1,\ldots,x_n)\in\mathbf{B}$
and $(x_1^n,y)\in P_a(g_1\op{i}g_2)$, where $y\in G$.

1) If $(x_1,\ldots,x_n)\in\mathbf{B_0}$, i.e.
$x_i=\mu_i^*(\opp{i_1}{i_s}y_1^s)$, $\,i=1,\ldots,n$, for some
$y_1,\ldots,y_s\in G\,$ and $\,i_1,\ldots,i_s\in\{1,\ldots,n\}$,
then evidently
\begin{equation}\label{21}
a\sqsubset y=(g_1\op{i\,}g_2)\opp{i_1}{i_s}y_1^s\, .
\end{equation}
But $i\in\{i,i_1,\ldots,i_s\}$, so
\begin{equation}\label{22}
\mu_i^*(\op{i\,}g_2\opp{i_1}{i_s}y_1^s)=\mu_i(\op{i\,}g_2\opp{i_1}{i_s}y_1^s)
=g_2\opp{i_1}{i_s}y_1^s\, .
\end{equation}
Therefore, by the $v$-negativity of $\chi$, from (\ref{21}) we
deduce
\[
a\sqsubset(g_1\op{i\,}g_2)\opp{i_1}{i_s}y_1^s\sqsubset
\mu_i(\op{i\,}g_2\opp{i_1}{i_s}y_1^s)=g_2\opp{i_1}{i_s}y_1^s\, .
\]
Hence $\,a\sqsubset g_2\opp{i_1}{i_s}y_1^s\, ,\,$ i.e.
$P_a(g_2)(x_1^n)=g_2\opp{i_1}{i_s}y_1^s$. For $k\neq i$ we have
obviously
\[
\mu_k^*(\op{i}g_2\opp{i_1}{i_s}y_1^s)=\mu_k^*(\opp{i_1}{i_s}y_1^s)=x_k,\quad
k=1,\ldots,n,
\]
which together with (\ref{22}) implies
$\,(x_1^{i-1},g_2\opp{i_1}{i_s}y_1^s,x_{i+1}^n)\in\mathbf{B_0}$.
Thus, by (\ref{21}), we obtain
\[
P_a(g_1)(x_1^{i-1},g_2\opp{i_1}{i_s}y_1^s,x_{i+1}^n)=
g_1\op{i\,}g_2\opp{i_1}{i_s}y_1^s=y,
\]
i.e. $P_a(g_1)\Big(x_1^{i-1},P_a(g_2)(x_1^n),x_{i+1}^n\Big)=y$.
So, $P_a(g_1)\op{i\,}P_a(g_2)(x_1^n)=y$, and in the consequence,
$(x_1^n,y)\in P_a(g_1)\op{i\,}P_a(g_2)$.

2) If $(x_1,\ldots,x_n)=(e_1,\ldots,e_n)$, then $a\sqsubset
y=g_1\op{i\,}g_2$. We have also $g_1\op{i\,}g_2\sqsubset g_2$,
$\,\mu_i^*(\op{i\,}g_2)=g_2$ and $\mu_k^*(\op{i\,}g_2)=e_k$ for
$k\neq i$, $\,k=1,\ldots,n$. From the above we obtain
$P_a(g_2)(e_1,\ldots,e_n)=g_2$,
$\,(e_1^{i-1},g_2,e_{i+1}^n)\in\mathbf{B_0}\,$ and
$P_a(g_1)(e_1^{i-1},g_2,e_{i+1}^n)=g_1\op{i\,}g_2$. Thus
$P_a(g_1)\Big(e_1^{i-1},P_a(g_2)(e_1^n),e_{i+1}^n\Big)=y$, whence
$(e_1^n,y)\in P_a(g_1)\op{i\,}P_a(g_2)$.

In this way, we have shown that in both cases $(x_1^n,y)\in
P_a(g_1)\op{i}P_a(g_2)$. This proves the inclusion
$\,P_a(g_1\op{i\,}g_2)\subset P_a(g_1)\op{i\,}P_a(g_2)$.

To prove the converse inclusion, let $(x_1^n,y)\in
P_a(g_1)\op{i\,}P_a(g_2)$. Then there exists $z\in G$ such that
\begin{eqnarray}
&\label{23} (x_1^n,z)\in P_a(g_2),&\\[4pt]
&\label{24} (x_1^{i-1}\,z\,x_{i+1}^n,y)\in P_a(g_1).&
\end{eqnarray}

1) If $(x_1,\ldots,x_n)\in\mathbf{B_0}$, then
$x_i=\mu_i^*(\opp{i_1}{i_s}y_1^s)$, $\,i=1,\ldots,n$, for some
$y_1,\ldots,y_s\in G$, $i_1,\ldots,i_s\in\{1,\ldots,n\}$. So, from
(\ref{23}) we get $a\sqsubset z=g_2\opp{i_1}{i_s}y_1^s$. Because
$\,\mu_i^*(\op{i\,}g_2\opp{i_1}{i_s}y_1^s)=g_2\opp{i_1}{i_s}y_1^s=z\,$
and
$\,\mu_k^*(\op{i\,}g_2\opp{i_1}{i_s}y_1^s)=\mu_k^*(\opp{i_1}{i_s}y_1^s)=x_k\,$
for $k\neq i$, $\,k=1,\ldots,n$, the condition (\ref{24}) can be
written in the form $a\sqsubset
y=g_1\op{i\,}g_2\opp{i_1}{i_s}y_1^s$, which is equivalent to
$(x_1^n,y)\in P_a(g_1\op{i\,}g_2)$.

2) If $(x_1,\ldots,x_n)=(e_1,\ldots,e_n)$, then (\ref{23}) gives
$a\sqsubset z=g_2$. Similarly, (\ref{24}) implies
$(e_1^{i-1},g_2,e_{i+1}^n)\in P_a(g_1)$. But
$(e_1^{i-1},g_2,e_{i+1}^n)\in\mathbf{B_0}$, therefore
$\,P_a(g_1)(e_1^{i-1},g_2,e_{i+1}^n)=g_1\op{i\,}g_2=y$, i.e.
$a\sqsubset y=g_1\op{i\,}g_2$, which means that $(e_1^n,y)\in
P_a(g_1\op{i}g_2)$.

So, in both cases we have $\, P_a(g_1)\op{i\,}P_a(g_2)\subset
P_a(g_1\op{i\,}g_2)$, which together with the previous inclusion
proves (\ref{20}). Thus, $P_a$ is a representation of a
$(2,n)$-semigroup $(G;\op{1\,},\ldots,\op{n})$ by $n$-place
functions.

Let $P_0$ be the sum of the family of representations $(P_a)_{a\in
G}$, i.e.
\begin{equation}\label{25}
  P_0=\sum\limits_{a\in G}P_a\, .
\end{equation}
Then, of course, $P_0$ also is a representation of this
$(2,n)$-semigroup by $n$-place functions.

Now we show that $\chi=\chi_{_{P_{0}}}$. Indeed, if
$(g_1,g_2)\in\chi_{_{P_{0}}}$, then, by (\ref{19}), we have
$(g_1,g_2)\in\bigcap\limits_{a\in G}\chi_{_{P_{a}}}\, ,\,$ i.e.
$\,\pr P_a(g_1)\subset\pr P_a(g_2)\,$ for every $\,a\in G$, which
means that
\[
(\forall a\in G)(\forall
x_1^n\in\mathbf{B})\Big((x_1,\ldots,x_n)\in\pr
P_a(g_1)\longrightarrow (x_1,\ldots,x_n)\in\pr P_a(g_2)\Big).
\]
This, in particular, for $(x_1,\ldots,x_n)=(e_1,\ldots,e_n)$ shows
that
\[
(\forall a\in G)\Big((\exists y\in G)(e_1^n,y)\in
P_a(g_1)\longrightarrow(\exists z\in G)(e_1^n,z)\in P_a(g_2)\Big),
\]
which means that
\[(\forall a\in G)(\forall y\in G)(a\sqsubset y=g_1\longrightarrow
(\exists z\in G)a\sqsubset z=g_2).\]
So,
\[
(\forall a\in G)(\forall g_{1}\in G)\Big(a\sqsubset g_1
\longrightarrow a\sqsubset g_2\Big).
\]
According to the reflexivity of $\chi$, the last condition implies
$\,g_1\sqsubset g_2$, i.e. $(g_1,g_2)\in\chi$. So,
$\chi_{_{P_0}}\subset\chi$.

Conversely, let $(g_1,g_2)\in\chi$, $a\in G$ and
$(x_1,\ldots,x_n)\in\pr P_a(g_1)$. If
$(x_1,\ldots,x_n)\in\mathbf{B_0}$, then $a\sqsubset
g_1\opp{i_1}{i_s}y_1^s$, where $x_i=\mu_i^*(\opp{i_1}{i_s}y_1^s)$,
$i=1,\ldots,n$. Since $\chi$ is $l$-regular, $g_1\sqsubset g_2$
implies $g_1\opp{i_1}{i_s}y_1^s\sqsubset g_2\opp{i_1}{i_s}y_1^s$.
Thus $a\sqsubset g_2\opp{i_1}{i_s}y_1^s$, i.e.
$(x_1,\ldots,x_n)\in\pr P_a(g_2)$. If
$(x_1,\ldots,x_n)=(e_1,\ldots,e_n)$, then $a\sqsubset g_1$.
Therefore $a\sqsubset g_2$, which gives $(e_1,\ldots,e_n)\in\pr
P_a(g_2)$. Thus, we have shown, that for any
$(x_1,\ldots,x_n)\in\mathbf{B}$ from $(x_1,\ldots,x_n)\in\pr
P_a(g_1)$ it follows $(x_1,\ldots,x_n)\in\pr P_a(g_2)$. From this,
according to (\ref{19}) and (\ref{25}), we conclude $\pr
P_0(g_1)\subset\pr P_0(g_2)$, i.e. $(g_1,g_2)\in\chi_{_{P_0}}$.
So, $\chi\subset\chi_{_{P_0}}$. Hence $\,\chi=\chi_{_{P_0}}$.

Since a $(2,n)$-semigroup $(G;\op{1\,},\ldots,\op{n})$ satisfies
the condition (\ref{14}), by Theorem~\ref{T2}, there exists an
isomorphic representation $\Lambda$ of this $(2,n)$-semigroup by
full $n$-place functions. Hence $\chi_{_{\Lambda}}=G\times G$ and
$\varepsilon_{_{\Lambda}}=\triangle_{G}$. Now consider the
representation $P$ of the given $(2,n)$-semigroup, which is
defined by the equality $\,P=\Lambda+P_0$. We have
$\chi_{_{P}}=\chi_{_{\Lambda}}\cap\chi_{_{P_0}}=G\times
G\cap\chi=\chi$ and $\varepsilon_{_{P}}=
\varepsilon_{_{\Lambda}}\cap\varepsilon_{_{P_0}}=
\bigtriangleup_G\cap\varepsilon_{_{P_0}}=\bigtriangleup_G$. This
means that $P$ is a faithful representation. So
$(G;\op{1\,},\ldots,\op{n},\chi)$ is isomorphic to some projection
quasi-ordered $(2,n)$-semigroup of $n$-place functions.
\end{proof}

\begin{minipage}{60mm}
\begin{flushleft}
{\it Dudek~W.~A.\\
 Institute of Mathematics\\
 Technical University\\50-370
Wroclaw \\ Poland\\
 Email: dudek@im.pwr.wroc.pl}
\end{flushleft}
\end{minipage}
\hfill
\begin{minipage}{60mm}
\begin{flushleft}
{\it Trokhimenko~V.~S.\\
 Department of Mathematics\\
Pedagogical University\\
21100 Vinnitsa \\
Ukraine\\
Email: vtrokhim@sovamua.com}
\end{flushleft}
\end{minipage}
\end{document}